\documentclass{article}
\usepackage{amssymb}

\usepackage{amsmath}

\newtheorem{theorem}{Theorem}

\newtheorem{corollary}{Corollary}

\newtheorem{definition}{Definition}
\newtheorem{example}[theorem]{Example}

\newtheorem{lemma}{Lemma}

\newtheorem{proposition}{Proposition}

\begin{document}

\title{Free point processes and free extreme values}
\author{G. Ben Arous \thanks{%
Courant Institute of Mathematical Sciences, 251 Mercer Street, New York, NY
10012, USA. e-mail: benarous@cims.nyu.edu. Research supported in part by NSF
Grants DMS 0808180, and OISE 0730136. } \ and V. Kargin \thanks{%
Department of Mathematics, Stanford University, Palo Alto, CA 94305, USA.
e-mail: kargin@stanford.edu}}
\maketitle

\begin{center}
\textbf{Abstract}
\end{center}

\begin{quotation}
We continue here the study of free extreme values begun in \cite%
{benarous_voiculescu06}. We study the convergence of the free point
processes associated with free extreme values to a free Poisson random
measure (\cite{voiculescu98}, \cite{bnt05}). We relate this convergence to
the free extremal laws introduced in \cite{benarous_voiculescu06} and give
the limit laws for free order statistics.
\end{quotation}

\section{Introduction}

In classical probability theory, the theory of extreme values for i.i.d.
random variables is elementary and well understood. Recently, a similar
theory has been introduced in the context of free probability theory, in
which the role of independent random variables is played by freely
independent operators in a Hilbert space (\cite{benarous_voiculescu06}). The
asymptotic behavior of the maximum of N free operators is given in \cite%
{benarous_voiculescu06}, where the maximum is taken for the spectral order
relation on operators (\cite{ando89}, \cite{olson71}). The theory emerging
is then parallel to the classical theory for maxima of i.i.d. random
variables. In this paper, we make the next step in developing this parallel
picture. We study the behavior of the full point process of normalized free
extreme values. We show that it converges to a free Poisson random measure
as soon as the normalized free maximum converges. One should notice that the
notion of \textquotedblleft free order statistics\textquotedblright\ is not
readily available. Indeed, the notion of a \textquotedblleft second largest
statistic\textquotedblright\ is not at all clear. This difficulty is
mirrored in the nature of the limiting object. Free Poisson random variables
are not discrete. We will see (in Theorem \ref{theorem_free_extremal_limits}%
) that our main convergence theorem (Theorem \ref%
{theorem_convergence_equivalence}) leads to results with no classical
analogs for order statistics.

The basic element in both classical and free theory of extremes is a
probability measure $\mu .$ In the classical case, we take a sequence of
i.i.d. random variables $X_{i},$ distributed according to $\mu ,$ and
introduce their order statistics, i.e., order them in non increasing order: 
\begin{equation*}
X^{(0)}\geq X^{(1)}\geq X^{(2)}\geq ...\geq X^{(n-1)},
\end{equation*}%
so that $X^{(0)}$ is the maximum of the n-sample, $X^{(1)}$ the second
largest value and so on. The basic question is to describe the asymptotic
behavior of the distribution of these order statistics once properly
normalized, when n tends to $\infty $.

Let $F_{n,k}$ denotes the distribution function of the normalized order
statistics $\frac{X^{(k)}-b_{n}}{a_{n}}$, for well chosen normalization
constants $a_{n}$ and $b_{n}$

\begin{equation*}
F_{n,k}(t)=P\left[ \frac{X^{(k)}-b_{n}}{a_{n}}\leq t\right] .
\end{equation*}%
The first question addresses the behavior of the maximum, i.e., the
asymptotic behavior of $F_{n,0}$. It was shown in the classical works by %
\cite{frechet27}, \cite{fisher_tippett28}, and \cite{gnedenko43} that there
are only three types of possible limit laws, to which $F_{n,0}$ can weakly
converge. These laws (Weibull, Frechet or Gumbel) are called ``extreme value
distributions'':%
\begin{equation*}
\begin{array}{ll}
\mathrm{Type\quad I:} & G\left( x\right) =\exp \left( -e^{-x}\right) ,\qquad
-\infty <x<\infty ; \\ 
&  \\ 
\mathrm{Type\quad II:} & G\left( x\right) =\left\{ 
\begin{array}{ll}
0, & x\leq 0, \\ 
\exp \left( -x^{-\alpha }\right) ,\quad \mathrm{for\quad some}\quad \alpha
>0, & x>0;%
\end{array}%
\right. \\ 
&  \\ 
\mathrm{Type\quad III:} & G\left( x\right) =\left\{ 
\begin{array}{ll}
\exp \left( -\left( -x\right) ^{\alpha }\right) ,\quad \mathrm{for\quad some}%
\quad \alpha >0, & x\leq 0, \\ 
1, & x>0.%
\end{array}%
\right.%
\end{array}%
\end{equation*}%
Moreover, the nature of the max-domain of attraction of these extreme value
distributions is well known as well as the possible choices for
normalization constants (\cite{leadbetter_lindgren_rootzen83}, \cite%
{resnick87}).

In the free probability context, a sequence of free self-adjoint operators $%
X_{i}$ is taken, such that each of $X_{i}$ has the spectral probability
distribution $\mu .$ In recent work \cite{benarous_voiculescu06}, a maximum
operation was defined which maps any $n$-tuple of self-adjoint operators to
another self-adjoint operator, which is called their maximum. The definition
is based on the so-called spectral order for self-adjoint operators: $%
A\preceq B$ \ iff all spectral projections $1_{\left( -\infty ,t\right]
}\left( A\right) $ are greater than or equal to the corresponding spectral
projections $1_{\left( -\infty ,t\right] }\left( B\right) .$

The spectral order is stronger than the usual order on operators, according
to which $A\leq B$ iff $B-A$ is non-negative definite. The main benefit of
the spectral order is that the set of all self-adjoint operators forms a
lattice with respect to this order. In particular, if $S$ is the set of all
operators $C$ such that $A_{i}\preceq C$ for each of $A_{1},\ldots ,A_{n}$,
then $S$ has a unique minimal element which is called $\max \left\{
A_{1},\ldots ,A_{n}\right\} .$ This property does not hold if self-adjoint
operators are considered with respect to the usual order on operators. Note,
however, that the lattice of selfadjoint operators with respect to the
spectral order is not a vector lattice in the sense that $A-B\succeq 0$ does
not imply that $A\succeq B.$ \ For a counter-example and other information
about the spectral order, see \cite{olson71}.

By analogy with the classical case, the sequence of normalized maxima is
defined as 
\begin{equation*}
\max_{1\leq i\leq n}\left\{ \left( X_{i}-b_{n}I\right) /a_{n}\right\}
\end{equation*}%
where the maximum here is understood with respect to the spectral order.
Then, $F_{n,0}^{free}\left( x\right) $ is defined as the spectral
distribution function of this normalized maximum.

In \cite{benarous_voiculescu06} the following question is solved: When does
the sequence of $F^{free}_{n,0} $ converges weakly?

The answer to this question is very similar to the answer in the classical
case: There are only three possible types of limit laws, and for a given $%
\mu ,$ the distributions $F_{n,0}^{free}$ can converge to only one of them:%
\begin{equation*}
\begin{array}{ll}
\mathrm{Type\quad I:} & G^{free}\left( x\right) =\left\{ 
\begin{array}{ll}
0, & x\leq 0, \\ 
1-e^{-x}, & x>0;%
\end{array}%
\right. \\ 
&  \\ 
\mathrm{Type\quad II:} & G^{free}\left( x\right) =\left\{ 
\begin{array}{ll}
0, & x\leq 1, \\ 
1-x^{-\alpha },\quad \mathrm{for\quad some}\quad \alpha >0, & x>1;%
\end{array}%
\right. \\ 
&  \\ 
\mathrm{Type\quad III:} & G^{free}\left( x\right) =\left\{ 
\begin{array}{ll}
0, & x\leq -1 \\ 
1-\left( -x\right) ^{\alpha },\quad \mathrm{for\quad some}\quad \alpha >0, & 
-1<x\leq 0, \\ 
1, & x>0.%
\end{array}%
\right.%
\end{array}%
\end{equation*}%
As in the classical case, this allows defining domains of attraction of the
free limit laws. Similar to the results about sums of free operators (\cite%
{bercovici_pata99}), an important fact is that, even though the limit laws
are different in the classical and free cases, the domains of attraction are
the same as well as the normalization constants! More precisely $F_{n,0}$
converges weakly to the extreme value distribution $G(x)$ iff $%
F_{n,0}^{free} $ converges weakly to $G^{free}$ of the same type as $G\left(
x\right) .$

This rigid link between classical and free probability theory for extreme
values is thus exactly similar to the analogous results for sums of i.i.d.
random variables, as developed in (\cite{bercovici_pata99}).

In order to investigate this situation further, let us return to the
classical case and consider the random point process 
\begin{equation*}
N_{n}=\sum_{i=1}^{n}\delta _{\left( X_{i}-b_{n}\right) /a_{n}}.
\end{equation*}%
The next question of classical extreme value theory is to understand the
convergence of this point process. This question is naturally related to the
convergence of the distributions $F_{n,k}.$ If $\mu $ is in the domain of
attraction of a classical extreme value distribution $G(x)$, or equivalently
if $F_{n,0}$ converges to $G(x)$ for some choice of normalization constants $%
a_{n}$ and $b_{n}$, then the point process $N_{n}$ weakly converges to a
Poisson random measure with intensity measure $\lambda \left( dx\right) $
with $\lambda (x,\infty )=-\log G(x)$. Conversely, if $N_{n}$ weakly
converges to a Poisson random measure with the intensity measure $\lambda
\left( dx\right) ,$ then the distribution of any order statistics $F_{n,k}$
converges to a limit law $G_{(k)}$ which is easily computable from $\lambda
(dx)$ or equivalently from $G(x)$, see below or (\cite{resnick87}).

What is the free analogue of the point process $N_{n}$? To motivate our
definition, note that we can think about $N_{n}$ as a linear functional on
the space of bounded measurable functions: $\left\langle
N_{n},f\right\rangle =:\sum_{i=1}^{n}f\left( \left( X_{i}-b_{n}\right)
/a_{n}\right) .$ This functional takes values in the space of bounded random
variables. We will define a free point process analogously. We begin with a
slightly greater generality and associate a free random process to any
triangular array of free random variables.

Let $\overline{\mathcal{A}}$ be the set of densely-defined closed operators
affiliated with a von Neumann algebra $\mathcal{A}$, and let $\mathcal{B}%
_{\infty }\left( \mathbb{R}\right) $ denote the set of all bounded, Borel
measurable functions $f:\mathbb{R}\rightarrow \mathbb{R}$.

\begin{definition}
Let $X_{i,n}\in \overline{\mathcal{A}},$ ($i=1,\ldots ,n;$ $n=1,\ldots $) be
a triangular array of freely independent self-adjoint variables. Then the
free point process $M_{n}$ associated with the array $X_{i,n}$ is the
sequence of $\mathcal{A}$-valued functionals on $\mathcal{B}_{\infty }\left( 
\mathbb{R}\right) $, defined by the following formula: 
\begin{equation*}
\left\langle M_{n},f\right\rangle :=\sum_{i=1}^{n}f\left( X_{i,n}\right) .
\end{equation*}
\end{definition}

The triangular array of free variables that we use in applications to free
extremes is, of course, $X_{i,n}=\left( X_{i}-b_{n}\right) /a_{n},$ where $%
X_{i}$ is a sequence of free self-adjoint variables.

We can also define the concept of weak convergence of a free point process
as a weak-$\ast $ convergence of the corresponding functionals. In the
classical case, after a suitable scaling, the point process $N_{n}$
converges to a Poisson random measure. It turns out that in the
non-commutative case the free point process converges to a \emph{free
Poisson random measure}, which was recently defined in \cite{voiculescu98}
and \cite{bnt05}. The following three theorems are the main results of our
paper.

\begin{theorem}
\label{theorem_convergence_equivalence}Let $G(x)$ be a classical extreme
value distribution, i.e. a Gumbel, Frechet or Weibull distribution. Let $%
\overline{x}=\inf \left\{ x\text{ }:\text{ }G\left( x\right) >0\right\} $
and define a measure $\lambda \left( dx\right) $ on $\left[ \overline{x}%
,\infty \right) $ by the equality $\lambda \left( \left( x,\infty \right)
\right) =-\log G\left( x\right) $. The following statements are equivalent:%
\newline
(i) $\mu $ belongs to the domain of attraction of the classical extremal
limit law $G(x)$, i.e., for some constants $a_{n}$ and $b_{n}$ the
distribution $F_{n,0}$ converges weakly to $G(x)$;\newline
(ii) $\mu $ belongs to the domain of attraction of the free extremal limit
law $G^{free}$, i.e., for some constants $a_{n}$ and $b_{n}$ the spectral
distribution of the normalized free maximum, $F_{n,0}^{free}$ converges
weakly to $G^{free}$;\newline
(iii) For some $a_{n}$ and $b_{n},$ the point process $N_{n}$ weakly
converges on $\left( \overline{x},\infty \right) $ to the Poisson random
measure with intensity $\lambda \left( dx\right) $;\newline
(iv) For some $a_{n}$ and $b_{n},$ the free point process $M_{n}$ weakly
converges on $\left( \overline{x},\infty \right) $ to the free Poisson
random measure with intensity $\lambda \left( dx\right) .$
\end{theorem}

In case one of the equivalent conditions in Theorem \ref%
{theorem_convergence_equivalence} is satisfied, then all the normalisation
constants $a_{n}$ and $b_{n}$ can be taken to be the same in all four
statements.

The equivalences of (i) and (iii) follows from the results in \cite%
{resnick87} (see, e.g., Section 4.2.2 on page 209), and the equivalence of
(i) and (ii) was proved in \cite{benarous_voiculescu06}. Thus, we only need
to prove the equivalence of (i) and (iv).

The equivalence of (i) and (iv) will be seen, in Section \ref%
{section_proof_of_Theorem_1}, as a consequence of the following more general
result about convergence of free point processes. Recall that a measure is
called Radon if $\mu \left( K\right) <\infty $ for every compact $K.$

\begin{theorem}
\label{theorem_weak_convergence_compact_cont}\textit{Let }$X_{i,n}$\textit{\
be a triangular array of free, self-adjoint random variables and let the
spectral probability measure of }$X_{i,n}$\textit{\ be }$\mu _{n}.$\textit{\
Let }$\lambda $\textit{\ be a Radon measure on }$D\subseteq \mathbb{R}.$%
\textit{\ The free point process }$M_{n}$\textit{\ associated with the array 
}$X_{i,n}$\textit{\ converges weakly on }$D$\textit{\ to a free Poisson
random measure }$M$\textit{\ with the intensity measure }$\lambda $ if and
only if\textit{\ }%
\begin{equation}
n\mu _{n}\left( A\right) \rightarrow \lambda \left( A\right)
\end{equation}%
\textit{for every Borel set }$A\subseteq D$\textit{. }
\end{theorem}

We now want to show what Theorem \ref{theorem_convergence_equivalence}
implies for free order statistics. We begin by recalling basic facts about
the classical theory of extreme values. If the measure $\mu $ is in the
domain of attraction of the extreme value distribution $G(x)$, then as
mentioned above, the convergence of the point process $N_{n}$ implies easily
the convergence of order statistics. Indeed with the notations introduced
above, it is easy to relate the distribution $F_{n,k}$ of the normalized $k$%
-th order statistics to the point process $N_{n}$, through the basic
identity: 
\begin{equation*}
F_{n,k}\left( t\right) =P\left[ {\frac{X^{(k)}-b_{n}}{a_{n}}\leq t}\right] =P%
\left[ N_{n}(t,\infty )\leq k\right] =E\left[ 1_{[0,k]}(\left\langle
N_{n},1_{(t,\infty )}\right\rangle )\right] .
\end{equation*}

This implies easily that the distribution $F_{n,k}$ of the properly
normalized order statistics weakly converges to the distribution 
\begin{equation*}
G_{(k)}(t)=\sum_{j=0}^{k}e^{-\lambda (t,\infty )}\frac{\lambda (t,\infty
)^{j}}{j!}.
\end{equation*}

We now want to see how this translates in the free context. More precisely,
let $X_{1},...,X_{n}$ be freely independent self-adjoint variables with
(possibly different) distribution functions $F_{i}$. Consider $M_{n}$ the
free point process associated with the sequence $X_{i}$ and let 
\begin{equation*}
Y_{n}\left( t\right) :=\left\langle M_{n},1_{\left( t,\infty \right)
}\right\rangle =\sum_{i=1}^{n}1_{\left( t,\infty \right) }\left(
X_{i}\right) .
\end{equation*}

\begin{definition}
For every real $k\geq 0,$ we say that $F_{n,k}^{free}\left( t\right) :=E%
\left[ 1_{\left[ 0,k\right] }\left( Y_{n}\left( t\right) \right) \right] $
is the \emph{distribution function of the }$k$\emph{-th\ order statistic} of
the sequence $X_{1},\ldots ,X_{n},$ and that it is the $k$\emph{-th order
free extremal convolution} of the spectral distribution functions $F_{i}.$
\end{definition}

Note that the definition is valid not only for all integer $k$ but also for
all non-negative real $k.$

One question that immediately arises is whether we can define an operator,
for which the distribution $F_{n,k}^{free}\left( t\right) $ would be a
spectral distribution function? The answer to this question is positive. The
condition $t^{\prime }\geq t$ $\ $implies that $Y_{n}\left( t^{\prime
}\right) \leq Y_{n}\left( t\right) $ and $1_{\left[ 0,k\right] }\left(
Y_{n}\left( t^{\prime }\right) \right) \geq 1_{\left[ 0,k\right] }\left(
Y_{n}\left( t\right) \right) .$ Therefore, as $t$ grows, the operators $1_{%
\left[ 0,k\right] }\left( Y_{n}\left( t\right) \right) $ form an increasing
family of projections and we can use this family to construct the required
operator by the spectral resolution theorem.

\begin{definition}
For every real $k\geq 0,$ let 
\begin{equation*}
Z^{\left( k\right) }=\int t\text{ }d1_{\left[ 0,k\right] }\left( Y_{n}\left(
t\right) \right) .
\end{equation*}%
We call $Z^{\left( k\right) }$ the $k$-th\emph{\ order statistic }of the
family\emph{\ }$X_{i}.$
\end{definition}

From the construction it is clear that $F_{n,k}\left( t\right) $ is the
spectral distribution function of the operator $Z^{\left( k\right) }.$

In complete analogy with the classical case the limits of these free
extremal convolutions can be computed using the limits of free point
measures. If $G\left( x\right) $ is one of the classical limit laws, then we
use $G^{\left( -1\right) }\left( x\right) $ to denote the functional inverse
of $G\left( x\right) .$ Let 
\begin{eqnarray*}
t_{-}\left( k\right) &=&G^{\left( -1\right) }\left( \exp \left[ -\left( 1+%
\sqrt{k}\right) ^{2}\right] \right) , \\
t_{0}\left( k\right) &=&G^{\left( -1\right) }\left( \frac{1}{e}\right) , \\
t_{+}\left( k\right) &=&G^{\left( -1\right) }\left( \exp \left[ -\left( 1-%
\sqrt{k}\right) ^{2}\right] \right) .
\end{eqnarray*}%
Let $\lambda \left( t\right) =$ $-\log G\left( t\right) $ and $p_{t}\left(
\xi \right) =\left( 2\pi \xi \right) ^{-1}\sqrt{4\xi -\left( 1-\lambda
\left( t\right) +\xi \right) ^{2}}.$

\begin{theorem}
\label{theorem_free_extremal_limits}Suppose that measure $\mu $ belongs to
the domain of attraction of a (classical) limit law $G\left( x\right) $ and $%
a_{n}$, $b_{n}$ are the corresponding norming constants. Assume that $X_{i}$
are free self-adjoint variables with the spectral probability measure $\mu $
and let $F_{n,k}^{free}\left( t\right) $ denote the distribution of the $k$%
-th order statistic of the family $\left( X_{i}-b_{n}\right) /a_{n},$ where $%
i=1,\ldots ,n.$ Then, as $n\rightarrow \infty ,$ the distribution $%
F_{n,k}^{free}\left( t\right) $ converges to a limit, $F_{\left( k\right)
}\left( t\right) ,$ which is given by the following formula: 
\begin{equation*}
F_{\left( k\right) }\left( t\right) =\left\{ 
\begin{array}{cc}
0, & \text{if }t<t_{-}, \\ 
\int_{\left( 1-\sqrt{\lambda \left( t\right) }\right) ^{2}}^{k}p_{t}\left(
\xi \right) d\xi , & \text{if }t\in \left[ t_{-},t_{0}\right] , \\ 
1-\lambda _{t}+\int_{\left( 1-\sqrt{\lambda \left( t\right) }\right)
^{2}}^{k}p_{t}\left( \xi \right) d\xi , & \text{if }\left( t_{0},t_{+}\right]
, \\ 
1-\lambda \left( t\right) 1_{\left[ 0,1\right) }\left( k\right) , & \text{if 
}t>t_{+}.%
\end{array}%
\right.
\end{equation*}
\end{theorem}

It turns out that in the particular case of the $0$-order free extremal
convolutions, their limits coincide with the limits discovered in \cite%
{benarous_voiculescu06} (see Definition 6.8 and Theorems 6.9 and 6.11):%
\begin{eqnarray*}
F_{\left( 0\right) }^{I}\left( t\right) &=&\left( 1-e^{-t}\right) 1_{\left(
0,\infty \right) }\left( t\right) ; \\
F_{\left( 0\right) }^{II}\left( t\right) &=&\left( 1-\frac{1}{t^{\alpha }}%
\right) 1_{\left( 1,\infty \right) }\left( t\right) ;\text{ and} \\
F_{\left( 0\right) }^{III}\left( t\right) &=&\left( 1-\left| t\right|
^{\alpha }\right) 1_{\left( -1,0\right) }\left( t\right) +1_{\left[ 0,\infty
\right) }\left( t\right) ,
\end{eqnarray*}%
where $\alpha $ is a positive parameter.

While we were mainly motivated by trying to extend the classical
probabilistic phenomena to the setting of free probability, it is worth
mentioning that the theory of free extreme values is directly related to
natural operations on random matrices (see the recent preprint \cite%
{benaych_cabanal08}). The results of this paper can easily be translated in
the context of \cite{benaych_cabanal08}.

The rest of the paper is organized as follows. Section \ref%
{section_Preliminaries} gives a brief introduction to free probability
theory. Section \ref{section_proof_of_Theorem_1} proves Theorem 1 using
Theorem 2. Section \ref{section_convergence_fPrm} details the definition of
the convergence of free point process and proves Theorem 2. And Section \ref%
{section_proof_of_Theorem_3} proves Theorem 3.

\section{Preliminaries\label{section_Preliminaries}}

\subsection{Free Independence}

\begin{definition}
A\emph{\ }$W^{\ast }$-\emph{probability space} is a pair $\left( \mathcal{A}%
,E\right) ,$ where $\mathcal{A}$ is a von Neumann algebra of bounded linear
operators acting on elements of a complex separable Hilbert space and $E$ is
a faithful normal trace that satisfies the condition $E(I)=1$. Operators
affiliated with algebra $\mathcal{A}$ are called \emph{non-commutative
random variables,} or simply \emph{random variables}, and the functional $E$
is called the \emph{expectation.}
\end{definition}

If $P\left( d\lambda \right) $ is the spectral resolution associated with a
normal operator $A,$ then we can define a measure $\mu \left( d\lambda
\right) =E\left( P\left( d\lambda \right) \right) .$ It is easy to check
that $\mu $ is a probability measure supported on the spectrum of $A$. We
call this measure, $\mu ,$ the \emph{spectral probability measure associated
with operator} $A$ \emph{and expectation} $E.$

The most important concept in free probability theory is that of free
independence of non-commutative random variables. Let a set of r.v. $%
A_{1},\ldots ,A_{n}$ be given. With each of them we can associate an algebra 
$\mathcal{A}_{i},$ which is generated by $A_{i}$ and $A_{i}^{\ast }$; that
is, it is the weak topology closure of all polynomials in variables $A_{i}$
and $A_{i}^{\ast }.$ Let $\overline{A}_{i}$ denote an arbitrary element of
algebra $\mathcal{A}_{i}.$

\begin{definition}
The algebras $\mathcal{A}_{1},\ldots ,\mathcal{A}_{n}$ (and variables $%
A_{1},\ldots ,A_{n}$ that generate them) are said to be \emph{freely
independent} or \emph{free}, if the following condition holds: 
\begin{equation*}
E\left( \overline{A}_{i(1)}\ldots \overline{A}_{i(m)}\right) =0,
\end{equation*}%
provided that $E\left( \overline{A}_{i(s)}\right) =0$ and $i(s+1)\neq i(s)$
for every $s$.
\end{definition}

For more information about non-commutative probability spaces and free
operators we refer the reader to Sections 2.2 - 2.5 in the book \cite%
{voiculescu_dykema_nica92} by Voiculescu, Dykema and Nica.

If $X$ and $Y$ are two free self-adjoint random variables with spectral
probabilities measures $\mu $ and $\nu $ respectively, then we denote the
spectral probability measure of $X+Y$ as $\mu \boxplus \nu ,$ and call it the%
\emph{\ free additive convolution} of $\mu $ and $\nu .$

\subsection{Free Poisson random variables}

Let $X$ be a self-adjoint operator that has the so-called \emph{free Poisson
distribution} with parameter (\textquotedblleft intensity\textquotedblright
) $\lambda .$ The continuous part of this distribution is supported on the
interval $\left[ \left( 1-\sqrt{\lambda }\right) ^{2},\left( 1+\sqrt{\lambda 
}\right) ^{2}\right] $ and the density is 
\begin{equation*}
p_{\lambda }\left( x\right) =\frac{\sqrt{4x-\left( 1-\lambda +x\right) ^{2}}%
}{2\pi x}.
\end{equation*}%
In addition, if $\lambda <1,$ then there is also an atom at zero with the
probability weight $1-\lambda .$ We call such an operator $X$ a
(non-commutative) \emph{Poisson random variable} with intensity $\lambda $
and size $1.$

The sum of two freely independent Poisson random variables of intensities $%
\lambda _{1}$ and $\lambda _{2}$ is again a Poisson random variable of
intensity $\lambda _{1}+\lambda _{2}$ (see, for example, a remark on page103
in \cite{hiai_petz00}).

If we scale a non-commutative Poisson random variable by $a,$ then we get a
variable, which we call a \emph{scaled (non-commutative) Poisson random
variable} of intensity $\lambda $ and size $a.$

Non-commutative Poisson random variables arise when we convolve a large
number, $N,$ of Bernoulli distributions that put probability $\lambda /N$ on 
$1$ and probability $1-\lambda /N$ on $0.$ The following result is
well-known, see \cite{hiai_petz00}, \cite{nica_speicher06},or \cite%
{voiculescu98}.

\begin{proposition}
\label{theorem_Poisson_limit}Suppose $\mu _{n},$ ($n=1,2,...)$ is a sequence
of Bernoulli distributions, such that $\mu _{n}\left( \left\{ 1\right\}
\right) \sim \lambda /n$ and $\mu _{n}\left( \left\{ 0\right\} \right)
=1-\mu _{n}\left( \left\{ 1\right\} \right) .$ Define $\nu _{n}$ as follows:%
\begin{equation*}
\nu _{n}=\underset{n\text{ times}}{\underbrace{\mu _{n}\boxplus ...\boxplus
\mu _{n}}}.
\end{equation*}%
Then $\nu _{n}$ weakly converges to the free Poisson distribution with
intensity $\lambda $ and size $1.$
\end{proposition}

\subsection{Free Poisson random measure}

\begin{definition}
Let $\left( \Theta ,\mathcal{B},\nu \right) $ be a measure space, and put 
\begin{equation*}
\mathcal{B}_{0}=\left\{ B\in \mathcal{B}:\nu \left( B\right) <\infty
\right\} .
\end{equation*}%
Let further $\left( \mathcal{A},E\right) $ be a $W^{\ast }$-probability
space, and let $\mathcal{A}_{+}$ denote the cone of positive operators in $%
\mathcal{A}.$ Then a \emph{free Poisson random measure} (\emph{fPrm}) on $%
\left( \Theta ,\mathcal{B},\nu \right) $ with values in $\left( \mathcal{A}%
,E\right) $ is a mapping $M:\mathcal{B}_{0}\rightarrow \mathcal{A}_{+},$
with the following properties: \newline
(i) For any set $B$ in $\mathcal{B}_{0}$, $M\left( B\right) $ is a free
Poisson variable with parameter $\nu \left( B\right) .$\newline
(ii) If $r\in \mathbb{N}$, and $B_{1},...,B_{r}\in \mathcal{B}_{0}$ are
disjoint, then $M\left( B_{1}\right) ,...,M\left( B_{r}\right) $ are free.%
\newline
(iii) If $r\in \mathbb{N}$, and $B_{1},...,B_{r}\in \mathcal{B}_{0}$ are
disjoint, then $M\left( \cup _{j=1}^{r}B_{j}\right) =\sum_{j=1}^{r}M\left(
B_{j}\right) .$
\end{definition}

The existence of a free Poisson measure for arbitrary spaces $\left( \Theta ,%
\mathcal{B},\nu \right) $ and $\left( \mathcal{A},E\right) $ was shown in %
\cite{voiculescu98}\ and a different proof was given in \cite{bnt05}.

Let $f$ be a real-valued simple function in $L^{1}\left( \Theta ,\mathcal{B}%
,\nu \right) ,$ i.e, suppose that it can be written as 
\begin{equation*}
f=\sum_{i=1}^{r}a_{i}1_{B_{i}},
\end{equation*}%
for a system of disjoint $B_{i}\in \mathcal{B}_{0}.$ Then we define the
integral of $f$ with respect to a Poisson random measure $M$ as follows:%
\begin{equation*}
\int_{\Theta }f\text{ }dM=\sum_{i=1}^{r}a_{i}M\left( B_{i}\right) .
\end{equation*}%
It is possible to check that this definition is consistent. Moreover, as it
is shown in \cite{bnt05}, this concept can be extended to a larger class of
functions:

\begin{proposition}
Let $f$ be a real-valued function in $L^{1}\left( \Theta ,\mathcal{B},\nu
\right) $ and suppose that $s_{n}$ is a sequence of real valued simple $%
\mathcal{B}$-measurable functions, satisfying the condition that there
exists a positive $\nu $-integrable function $h\left( \theta \right) ,$ such
that $\left\vert s_{n}\left( \theta \right) \right\vert \leq h\left( \theta
\right) $ for all $n$ and $\theta .$ Suppose also that $\lim_{n\rightarrow
\infty }s_{n}\left( \theta \right) =f\left( \theta \right) $ for all $\theta 
$ $\ $Then integrals $\int_{\Theta }s_{n}$ $dM$ are well-defined and
converge in probability to a self-adjoint (possibly unbounded) operator $%
I\left( f\right) $ affiliated with $\mathcal{A}$. Furthermore, the limit $%
I\left( f\right) $ is independent of the choice of approximating sequence $%
s_{n}$ of simple functions$.$
\end{proposition}

The resulting functional $I\left( f\right) $ is defined for all real valued
functions $f$ in $L^{1}\left( \Theta ,\mathcal{B},\nu \right) $ and is
called the \emph{integral with respect to the free Poisson random measure} $%
M $. It possesses all the usual properties of the integral: additivity,
linear scaling, continuity, etc.

\section{Proof of Theorem 1 \label{section_proof_of_Theorem_1}}

As was noted in Introduction, only the equivalence of (i) and (iv) needs a
proof. The equivalence of (i) and (iv) can be reduced to a problem about
convergence of free point processes. Indeed, let $\mu _{n}\left( A\right)
=\mu \left( a_{n}A+b_{n}\right) .$ Then (i) is equivalent to the statement
that $n\mu _{n}\left( A\right) \rightarrow \lambda \left( A\right) $ for all
Borel sets $A\subset \left( \overline{x},\infty \right) .$

Indeed, suppose that $\mu $ is in the domain of attaction of $G\left(
x\right) ,$ and let $F\left( x\right) $ denote the distribution function of
the measure $\mu .$ Then

\begin{equation*}
F^{n}\left( a_{n}x+b_{n}\right) \rightarrow G\left( x\right) ,
\end{equation*}%
For every $x\in \left( \overline{x},\infty \right) ,$ $G\left( x\right) $ is
positive, hence we can take logarithms and get%
\begin{equation*}
n\log F\left( a_{n}x+b_{n}\right) \rightarrow \log G\left( x\right) ,
\end{equation*}%
which is equivalent to 
\begin{equation*}
n\left( 1-F\left( a_{n}x+b_{n}\right) \right) \rightarrow -\log G\left(
x\right) \equiv \lambda \left( \left( x,\infty \right) \right) .
\end{equation*}%
Consequently, 
\begin{equation*}
n\mu _{n}\left( \left( x,\infty \right) \right) \rightarrow \lambda \left(
\left( x,\infty \right) \right) ,
\end{equation*}%
from which we conclude that $n\mu _{n}\left( A\right) \rightarrow \lambda
\left( A\right) $ for all Borel sets $A\subset \left( \overline{x},\infty
\right) .$

By reversing the steps of this argument we obtain the reverse implication:
If $n\mu _{n}\left( A\right) \rightarrow \lambda \left( A\right) $ for all
Borel sets $A\subset \left( \overline{x},\infty \right) ,$ then $\mu $ is in
the domain of attaction of $G\left( x\right) $, and (i) holds.

Therefore the equivalence of (i) and (iv) follows from Theorem \ref%
{theorem_weak_convergence_compact_cont} if we take $\left(
X_{i}-a_{n}\right) /b_{n}$ as the triangular array $X_{i,n}.$

\section{Proof of Theorem 2 \label{section_convergence_fPrm}}

\subsection{Weak Convergence}

In this section, we define precisely the mode of convergence of free point
measures that we use. It corresponds to the weak convergence of point
processes in the classical case.

Let $D$ be a Borel subset of $\mathbb{R}$ and let $\mathcal{F}_{K}^{\infty
}\left( D\right) $ denote the space of bounded, Borel measurable functions
that have compact support on $D.$

\begin{definition}
We say that a free point process $M_{n}$ \emph{converges weakly on D} to a
free Poisson random measure $M,$ which is defined on $\left( D,\mathcal{B}%
,\lambda \right) $ and takes values in $\mathcal{A}$, if for every function $%
f\in \mathcal{F}_{K}^{\infty }\left( D\right) $ the following convergence
holds: 
\begin{equation*}
\left\langle M_{n},f\right\rangle \overset{d}{\rightarrow }\int_{\mathbb{R}}f%
\text{ }dM.
\end{equation*}
\end{definition}

Sometimes we also need to speak about convergence with respect to a class of
functions, which is different from $\mathcal{F}_{K}^{\infty }\left( D\right) 
$.

\begin{definition}
We say that a free point process $M_{n}$ \emph{converges weakly} \emph{with
respect to} \emph{a class of functions} $\mathcal{F}$ to a free Poisson
random measure $M,$ if for every function $f\in \mathcal{F}$ the following
convergence holds: 
\begin{equation*}
\left\langle M_{n},f\right\rangle \overset{d}{\rightarrow }\int_{\mathbb{R}}f%
\text{ }dM.
\end{equation*}
\end{definition}

We will prove Theorem \ref{theorem_weak_convergence_compact_cont} by
considering initially the convergence of free point processes $M_{n}$ with
respect to the class of simple functions (i.e., finite sums of indicator
functions), and then approximating functions from a more general class by
simple functions.

\subsection{Convergence with respect to simple functions}

Let $\mathcal{S}\left( D\right) $ be the class of simple functions on $%
D\subset \mathbb{R}$, i.e., the class of finite sums of indicator functions
of Borel sets belonging to $D$.

\begin{proposition}
\label{theorem_weak_convergence_simple} \textit{Let }$X_{i,n}$\textit{\ be a
triangular array of free, self-adjoint random variables and let the spectral
probability measure of }$X_{i,n}$\textit{\ be }$\mu _{n}.$\textit{\ Let }$%
\lambda $\textit{\ be a Radon measure on }$D\subseteq \mathbb{R}.$\textit{\
If}%
\begin{equation*}
n\mu _{n}\left( A\right) \rightarrow \lambda \left( A\right)
\end{equation*}%
for each Borel set $A\subset D,$ then the free point process $M_{n}$
associated with the array $X_{i,n}$ converges weakly with respect to $%
\mathcal{S}\left( D\right) $ to a free Poisson random measure $M$ with the
intensity measure $\lambda $.
\end{proposition}

Before proving this proposition, we derive some auxiliary results.

\begin{lemma}
\label{lemma_weak_convergence_characteristic} Suppose $X_{i,n}$ is an array
of free and identically distributed random variables with the spectral
measure $\mu _{n}.$ \ Let $n\mu _{n}\left( A\right) \rightarrow \lambda
(A)<\infty $ as $n\rightarrow \infty .$ Let $Z_{i,n}=1_{A}\left(
X_{i,n}\right) .$ Then as $n\rightarrow \infty ,$ the sum $%
S_{n}=\sum_{i=1}^{n}Z_{i,n}$ converges in distribution to a free Poisson
random variable with intensity $\lambda \left( A\right) .$
\end{lemma}

\textbf{Proof:} Note that $Z_{i,n}$ are projections with expectation $\mu
_{n}\left( A\right) $ and they are free. Therefore, $\sum_{i=1}^{n}Z_{i,n}$
is the sum of free projections and we can use Proposition \ref%
{theorem_Poisson_limit} to infer the claim of the lemma. QED.

As the next step to the proof of Proposition \ref%
{theorem_weak_convergence_simple} we need to check that if Borel sets $A_{k}$
are disjoint, then the sums $S_{k}=\sum_{i=1}^{n}1_{A_{k}}\left(
X_{i,n}\right) $ are asymptotically free with respect to growing $n.$

Recall the definition of the asymptotic freeness: Let $\left( \mathcal{A}%
_{i},E_{i}\right) $ be a sequence of non-commutative probability spaces and
let $X_{i}$ and $Y_{i}$ be two random variables in $\mathcal{A}_{i}.$ Let
also $x$ and $y$ be two free operators in a non-commutative probability
space $\left( \mathcal{A},E\right) $.

\begin{definition}
\label{definition_asymptotic_freeness}The sequences $X_{i}$ and $Y_{i}$ are
called \emph{asymptotically free} if the sequence of pairs $\left(
X_{i},Y_{i}\right) $ converges in distribution to the pair $\left(
x,y\right) .$ That is, for every $\varepsilon >0$ and every sequence of $k$%
-tuples $\left( n_{1},...,n_{k}\right) $ with non-negative integers $n_{j}$,
there exists such $i_{0}$ that for $i\geq i_{0},$ the following inequality
holds:%
\begin{equation*}
\left\vert E_{i}\left(
X_{i}^{n_{1}}Y_{i}^{n_{2}}...X_{i}^{n_{k-1}}Y_{i}^{n_{k}}\right) -E\left(
x^{n_{1}}y^{n_{2}}...x^{n_{k-1}}y^{n_{k}}\right) \right\vert \leq
\varepsilon .
\end{equation*}
\end{definition}

At the cost of more complicated notation, this definiton can be generalized
to the case of more than two variables.

\begin{lemma}
\label{theorem_sums_asymptotic_freeness copy(1)}Let $P_{i,n}^{\left(
k\right) },$ (where $n=1,2,...;$ $i=1,...,n,$ and $k=1,...,r$) be
projections of dimension $\lambda ^{\left( k\right) }/n.$ Assume that for
each $n,$ algebras $\mathcal{A}_{i}$ generated by sets $\left\{
P_{i,n}^{\left( k\right) }\right\} _{k=1}^{r}$ are free. Also assume that
for each $n$ and $i,$ the projections $P_{i,n}^{\left( k\right) }$ are
orthogonal to each other, i.e., $P_{i,n}^{\left( k\right) }P_{i,n}^{\left(
k^{\prime }\right) }=0$ for every pair $k\neq k^{\prime }.$ Let $%
S_{n}^{\left( k\right) }=\sum_{i=1}^{n}P_{i,n}^{\left( k\right) }.$ Then as $%
n\rightarrow \infty ,$ the sequences $S_{n}^{\left( k\right) }$ converge in
distribution to freely independent variables $S^{\left( k\right) }$ that
have free Poisson distributions with parameters $\lambda ^{\left( k\right)
}, $ respectively. In particular, the sequences $S_{n}^{\left( k\right) }$
are asymptotically free with respect to growth in $n$.
\end{lemma}

\textbf{Proof: }The fact that each of the sequences $S_{n}^{\left( k\right)
} $ converge in distribution to a variable $S^{\left( k\right) }$ that has a
free Poisson distribution is clear from Proposition \ref%
{theorem_Poisson_limit}. The essential part is to prove that asymptotic
freeness holds. This claim is a direct consequence of Speicher's
multidimensional limit theorem (see, for example, Theorem 13.1 in the book %
\cite{nica_speicher06} by Nica and Speicher). Indeed, we need to prove that
all mixed free cumulants of the limit are zero. By Speicher's theorem, this
is equivalent to the statement that the following limits are zero:%
\begin{equation*}
\lim_{n\rightarrow \infty }nE\left( P_{1,n}^{\left( k_{1}\right)
}P_{1,n}^{\left( k_{2}\right) }\ldots P_{1,n}^{\left( k_{s}\right) }\right)
=0.
\end{equation*}%
Here $k_{1},\ldots ,k_{s}$ is an arbitrary $s$-tuple with the property that
it has a pair of distinct coordinates, i.e., $k_{i}\neq k_{j}.$ However, the
fact that these limits are zero is clear from the assumption that the
projections $P_{i,n}^{\left( k\right) }$ are orthogonal to each other. QED.

Now we can proceed to the proof of Proposition \ref%
{theorem_weak_convergence_simple}.

\textbf{Proof: }Let $f=\sum_{k=1}^{r}c_{k}1_{A_{k}}\left( x\right) ,$ where $%
A_{k}$ are disjoint Borel sets. Using the assumption that $n\mu _{n}\left(
A_{k}\right) \rightarrow \lambda \left( A_{k}\right) $ and Lemma \ref%
{lemma_weak_convergence_characteristic}, we can find a free Poisson random
measure $M$ such that%
\begin{equation*}
\sum_{i=1}^{n}1_{A_{k}}\left( X_{i,n}\right) \overset{d}{\rightarrow }%
M\left( A_{k}\right) =\int_{\mathbb{R}}1_{A_{k}}\left( x\right) M\left(
dx\right)
\end{equation*}%
as $n\rightarrow \infty .$ Indeed, it is enough to take a Poisson random
measure $M$ with the intensity measure $\lambda .$

In addition, by Lemma \ref{theorem_sums_asymptotic_freeness copy(1)}, sums $%
S_{k}=\sum_{i=1}^{n}1_{A_{k}}\left( X_{i,n}\right) $ become asymptotically
free for different $k$ as $n$ grows. Since $M\left( A_{k}\right) $ are free
by the definition of the free Poisson measure, this implies that 
\begin{equation*}
\sum_{k=1}^{r}c_{k}\sum_{i=1}^{n}1_{A_{k}}\left( X_{i,n}\right) \overset{d}{%
\rightarrow }\sum_{k=1}^{r}c_{k}M\left( A_{k}\right)
=\sum_{k=1}^{r}c_{k}\int_{\mathbb{R}}1_{A_{k}}\left( x\right) M\left(
dx\right) .
\end{equation*}%
as $n\rightarrow \infty .$ Therefore, 
\begin{equation*}
\sum_{i=1}^{n}f\left( X_{i,n}\right) \overset{d}{\rightarrow }\int_{\mathbb{R%
}}f\left( x\right) M\left( dx\right) ,
\end{equation*}%
where we used the additivity property of the integral with respect to a free
Poisson random measure (see \cite{bnt05}, Remark 4.2(b)). QED.

\subsection{Convergence with respect to bounded, Borel measurable functions
with compact support}

The goal of this section is to prove our main Theorem \ref%
{theorem_weak_convergence_compact_cont}.

Consider a bounded, Borel measurable, compactly supported function $%
f:D\rightarrow \mathbb{R},$ such that $0\leq f\leq 1.$ (A more general case
of a function $f,$ which satisfies $C_{1}\leq f\leq C_{2},$ can be treated
similarly.) For positive integers $N=1,2,\ldots ,$ and $k=1,\ldots ,N,$
define the set 
\begin{equation*}
A_{k}^{\left( N\right) }=\left\{ x\in \mathrm{supp}\left( f\right) :\frac{k-1%
}{N}<f\left( x\right) \leq \frac{k}{N}\right\} .
\end{equation*}%
The sets $A_{k}^{\left( N\right) }$ are disjoint, measurable, and have
finite $\lambda $-measure. Their union is $D.$

We define lower and upper approximations to the function $f$ as follows:%
\begin{equation*}
l^{N}\left( x\right) =\sum_{k=1}^{N}\frac{k-1}{N}1_{A_{k}^{\left( N\right)
}}\left( x\right) ,
\end{equation*}%
and 
\begin{equation*}
u^{N}\left( x\right) =\sum_{k=1}^{N}\frac{k}{N}1_{A_{k}^{\left( N\right)
}}\left( x\right) ,
\end{equation*}

We note that:\newline
(i) $l^{N}\left( x\right) \leq u^{N}\left( x\right) ;$\newline
(ii) $l^{N}\left( x\right) $ is an increasing sequence of functions;\newline
(iii) $u^{N}\left( x\right) $ is a decreasing sequence of functions, and%
\newline
iv) $\lim_{N\rightarrow \infty }l^{N}\left( x\right) -u^{N}\left( x\right)
=0 $ uniformly in $x.$

The functions $l^{N}\left( x\right) $ and $u^{N}\left( x\right) $ are
simple: $l^{N}\left( x\right) =\sum_{i=1}^{N}c_{k}^{\left( N\right)
}1_{A_{k}^{\left( N\right) }}\left( x\right) $ and $u^{N}\left( x\right)
=\sum_{i=1}^{N}d_{k}^{\left( N\right) }1_{A_{k}^{\left( N\right) }}\left(
x\right) .$ Note also that $\sup_{k}\left( d_{k}^{\left( N\right)
}-c_{k}^{\left( N\right) }\right) =1/N$ converges to zero as $N\rightarrow
\infty .$

Let us drop for convenience the superscript $N$ when we consider it as
fixed, and simply write $l\left( x\right)
=\sum_{i=1}^{N}c_{k}1_{A_{k}}\left( x\right) $ and $u\left( x\right)
=\sum_{i=1}^{N}d_{k}1_{A_{k}}\left( x\right) ,$ where $A_{k}$ are disjoint
Borel-measurable sets. By Proposition \ref{theorem_weak_convergence_simple},
as $n\rightarrow \infty ,$ 
\begin{equation*}
\sum_{i=1}^{n}l\left( X_{i,n}\right) \overset{d}{\rightarrow }%
\sum_{k=1}^{N}c_{k}M_{k},
\end{equation*}%
where $M_{k}$ are freely independent Poisson random variables with
intensities $\lambda _{k}=\lambda \left( A_{k}\right) $. Let $F_{l}\left(
x\right) $ denote the distribution function of $\sum_{k=1}^{N}c_{k}M_{k}.$

Similarly, 
\begin{equation*}
\sum_{i=1}^{n}u\left( X_{i,n}\right) \overset{d}{\rightarrow }%
\sum_{k=1}^{N}d_{k}M_{k},
\end{equation*}%
and we denote the distribution function of $\sum_{k=1}^{N}d_{k}M_{k}$ as $%
F_{u}\left( x\right) $.

Let $F_{f,n}$ denote the distribution function of $\sum_{i=1}^{n}f\left(
X_{i,n}\right) $ and let $F_{f}$ be one of the limit points of this sequence
of distribution functions.

\begin{proposition}
\label{proposition_limit_point_is_between}$F_{f}$ is a distribution function
and $F_{u}\left( x\right) \leq F_{f}\left( x\right) \leq F_{l}\left(
x\right) $ for every $x.$
\end{proposition}

\textbf{Proof:} We will infer this from Lemma \ref%
{lemma_Weyl_inequalities_II} below and its Corollary. This lemma is a
particular case of Weyl's eigenvalue inequalities for operators in a von
Neumann algebra of type $II_{1}.$ If $F_{A}\left( x\right) $ is the spectral
distribution function of a self-adjoint operator $A,$ then we define the
eigenvalue function $\theta _{A}\left( t\right) =\inf \left\{ x:F_{A}\left(
x\right) \geq 1-t\right\} .$ The function $\theta _{A}\left( t\right) $ is
non-increasing and right-continuous. Intuitively, it can be thought of as a
``sequence of eigenvalues'' of $A,$ indexed in decreasing order by parameter 
$t.$

Let us use notation $\theta _{A}\left( t-0\right) $ to denote $%
\lim_{\varepsilon \downarrow 0}\theta _{A}\left( t-\varepsilon \right) .$
Then the following generalization of Weyl inequalities holds:

\begin{lemma}
\label{lemma_Weyl_inequalities_II}If $A$ and $B$ are two bounded
self-adjoint operators from a $W^{\ast }$-probability space $\mathcal{A}$
and if $B$ is non-negative definite, then%
\begin{eqnarray*}
\theta _{A}\left( t\right) &\leq &\theta _{A+B}\left( t\right) \leq \theta
_{A}\left( t\right) +\left\| B\right\| ,\text{ and} \\
\theta _{A}\left( t-0\right) &\leq &\theta _{A+B}\left( t-0\right) \leq
\theta _{A}\left( t-0\right) +\left\| B\right\| .
\end{eqnarray*}
\end{lemma}

\begin{corollary}
\label{corollary_domination}If $B\geq 0,$ then $\mu _{A+B}\gg \mu _{A},$
that is, $F_{A+B}\left( x\right) \leq F_{A}\left( x\right) $ for each $x.$
\end{corollary}

\textbf{Proof of Lemma \ref{lemma_Weyl_inequalities_II}:} These results
easily follow from an inequality in \cite{bercovici_li01} which states that
if $\left( a-\varepsilon ,a\right) \subset \left[ 0,1\right] $, $\left(
b-\varepsilon ,b\right) \subset \left[ 0,1\right] ,$ and $a+b\leq 1$, then%
\begin{equation}
\int_{a+b-\varepsilon }^{a+b}\theta _{A+B}\left( t\right) dt\leq
\int_{a-\varepsilon }^{a}\theta _{A}\left( t\right) dt+\int_{b-\varepsilon
}^{b}\theta _{B}\left( t\right) dt.  \label{formula_weyl_bercovici_li}
\end{equation}%
QED.

By Corollary \ref{corollary_domination}, for each $n$ the distribution $%
F_{f,n}$ is between the distribution functions of $\sum_{i=1}^{n}u\left(
X_{i,n}\right) $ and $\sum_{i=1}^{n}l\left( X_{i,n}\right) .$ As $n$ grows,
these two sequences of distribution functions approach $F_{u}\left( x\right) 
$ and $F_{l}\left( x\right) ,$ respectively. Therefore, every limit point of 
$F_{f,n}$ is between $F_{u}$ and $F_{l}.$ The claim that $F_{f}$ is a
distribution function follows from the fact that both $F_{u}$ and $F_{l}$
are distribution functions. QED.

Now we want to show that $F_{u}^{\left( N\right) }\left( x\right) $
approaches $F_{l}^{\left( N\right) }\left( x\right) $ as $N$ grows.

Recall that the\emph{\ Levy distance} between two distribution functions is
defined as follows:%
\begin{equation*}
d_{L}\left( F_{A},F_{B}\right) =\sup_{x}\inf \left\{ s\geq 0:F_{B}\left(
x-s\right) -s\leq F_{A}\left( x\right) \leq F_{B}\left( x+s\right) +s\text{ }%
\right\} .
\end{equation*}

We can interpret this distance geometrically. Let $\Gamma _{A}$ be the graph
of function $F_{A},$ and at the points of discontinuity let us connect the
left and right limits by a (vertical) straight line interval. Call the
resulting curve $\widetilde{\Gamma }_{A}.$ Similarly define $\widetilde{%
\Gamma }_{B}.$ Let $d$ be the maximum distance between $\widetilde{\Gamma }%
_{A}$ and $\widetilde{\Gamma }_{B}$ in the direction from the south-east to
the north-west, i.e., in the direction which is obtained by rotating the
vertical direction by $\pi /4$ counter-clockwise. Then $d_{L}\left(
F_{A},F_{B}\right) =d/\sqrt{2}.$

\begin{proposition}
\label{lemma_Levy_distance_difference_coefficients} Let $K$ be the sum of
intensities of freely independent Poisson random variables $M_{k}$ and let $%
F_{l}\left( x\right) $ and $F_{u}\left( x\right) $ be distribution functions
of $\sum_{k=1}^{N}c_{k}M_{k}$ and $\sum_{k=1}^{N}d_{k}M_{k}$ Then%
\begin{equation*}
d_{L}\left( F_{l},F_{u}\right) \leq \left( 2K+3\sqrt{K}+1\right) \sup_{1\leq
k\leq N}\left( d_{k}-c_{k}\right) .
\end{equation*}
\end{proposition}

\textbf{Remark:} In the proof of Theorem \ref%
{theorem_weak_convergence_compact_cont}, the finiteness of $K$ will be
ensured by the assumptions that measure $\lambda $ is Radon and that $f$ has
a compact support.

For the proof of this proposition we need two lemmas. Lemma \ref%
{lemma_bound_lin_comb_Poisson_rv} provides a bound on the norm of the sum of
scaled Poisson random variables in terms of the sizes of these variables,
and Lemma \ref{lemma_Levy_distance_norm_difference} relates the Levy
distance between two random variables to the norm of their difference.

\begin{lemma}
\label{lemma_bound_lin_comb_Poisson_rv}Let $M_{i},$ ($i=1,...,r$) be freely
independent Poisson random variables, which have intensities $\lambda _{i},$
and let $b_{i}$ be non-negative real numbers. Assume that $%
\sum_{i=1}^{r}\lambda _{i}\leq K$ and let $b=\sup_{1\leq i\leq r}b_{i}.$
Then 
\begin{equation*}
\left\Vert \sum_{i=1}^{r}b_{i}M_{i}\right\Vert \leq b\left( 2K+3\sqrt{K}%
+1\right) .
\end{equation*}
\end{lemma}

\textbf{Proof:} Let $X_{i}$ be free self-adjoint random variables that have
zero mean. Then by an inequality from \cite{voiculescu86}:%
\begin{equation*}
\left\Vert \sum_{i=1}^{r}X_{i}\right\Vert \leq \max_{1\leq i\leq
r}\left\Vert X_{i}\right\Vert +\sqrt{\sum_{i=1}^{r}Var\left( X_{i}\right) }.
\end{equation*}%
If $Y_{i}$ are free self-adjoint random variables with non-zero mean, and $%
X_{i}=Y_{i}-E\left( Y_{i}\right) ,$ then the previous inequality implies
that 
\begin{eqnarray}
\left\Vert \sum_{i=1}^{r}Y_{i}\right\Vert &\leq &\left\vert
\sum_{i=1}^{r}E\left( Y_{i}\right) \right\vert +\left\Vert
\sum_{i=1}^{r}X_{i}\right\Vert  \notag \\
&\leq &\left\vert \sum_{i=1}^{r}E\left( Y_{i}\right) \right\vert
+\max_{1\leq i\leq r}d\left( Y_{i}\right) +\sqrt{\sum_{i=1}^{r}Var\left(
Y_{i}\right) },  \label{inequality_sum_Poisson_rv}
\end{eqnarray}%
where $d\left( Y_{i}\right) $ is the diameter of the support of $Y_{i}.$

We will apply this inequality to $Y_{i}=b_{i}M_{i}$ and estimate each of the
three terms on the right-hand side of (\ref{inequality_sum_Poisson_rv}) in
turn:

1) Since $E\left( M_{i}\right) =\lambda _{i},$ and $\sum \lambda _{i}\leq K,$
therefore $\sum_{i=1}^{r}b_{i}E\left( M_{i}\right) \leq bK.$

2) The diameter of the support of $b_{i}M_{i}$ is less or equal to $%
b_{i}\left( 1+\sqrt{\lambda _{i}}\right) ^{2}\leq b\left( 1+2\sqrt{K}%
+K\right) .$

3) Since $Var\left( M_{i}\right) =\lambda _{i},$ therefore $\sqrt{%
\sum_{i=1}^{r}Var\left( b_{i}M_{i}\right) }\leq b\sqrt{K}.$

In sum, $\left\Vert \sum_{i=1}^{r}b_{i}M_{i}\right\Vert \leq b\left( 2K+3%
\sqrt{K}+1\right) .$ QED.

\begin{lemma}
\label{lemma_Levy_distance_norm_difference}Let $A$ and $B$ be two bounded
self-adjoint operators from a $W^{\ast }$-probability space $\mathcal{A}$
and assume that $B-A\geq 0$. Then%
\begin{equation*}
d_{L}\left( F_{A},F_{B}\right) \leq \left\Vert B-A\right\Vert .
\end{equation*}
\end{lemma}

\textbf{Proof:} Let $F_{A}$ and $F_{B}$ be distribution functions, and $%
\theta _{A}$ and $\theta _{B}$ be the corresponding eigenvalue functions.
Then we claim that%
\begin{equation}
d_{L}\left( F_{A},F_{B}\right) \leq \sup_{0\leq t\leq 1}\left| \theta
_{A}\left( t\right) -\theta _{B}\left( t\right) \right| .
\label{inequality_Levy_lambda}
\end{equation}%
Indeed, let the graphs of functions $\theta _{A}$ and $\theta _{B}$ be
denoted as $\Lambda _{A}$ and $\Lambda _{B},$ respectively. Connecting the
left and right limits at the points of discontinuity gives us the curves $%
\widetilde{\Lambda }_{A}$ and $\widetilde{\Lambda }_{B}.$ It is easy to see
that these curves can be obtained from curves $\widetilde{\Gamma }_{A}$ and $%
\widetilde{\Gamma }_{B}$ (i.e., the graphs of $F_{A}\left( x\right) $ and $%
F_{B}\left( x\right) $ with connected limits at the points of discontinuity)
by rotating them around the point $\left( 0,1\right) $ counter-clockwise by
the angle $\pi /2$ and then shifting the result of the rotation by vector $%
\left( 0,-1\right) .$ It follows that the distance $d,$ which was used in
the definition of the Levy distance can also be defined as the maximum
distance between $\widetilde{\Lambda }_{A}$ and $\widetilde{\Lambda }_{B}$
in the direction from the south-west to the north-east, i.e., in the
direction which is obtained by rotating the vertical direction by $\pi /4$
clockwise.

Since $\theta _{A}\left( t\right) $ and $\theta _{B}\left( t\right) $ are
non-increasing functions, therefore 
\begin{equation*}
d\leq \sqrt{2}\sup_{0\leq t\leq 1}\left| \theta _{A}\left( t\right) -\theta
_{B}\left( t\right) \right| .
\end{equation*}%
This implies $d_{L}\left( F_{A},F_{B}\right) \leq \sup_{0\leq t\leq 1}\left|
\theta _{A}\left( t\right) -\theta _{B}\left( t\right) \right| .$

Inequality (\ref{inequality_Levy_lambda}) and Lemma \ref%
{lemma_Weyl_inequalities_II} imply the statement of the lemma. QED.

Now we can prove Proposition \ref%
{lemma_Levy_distance_difference_coefficients}:

\textbf{Proof of Proposition \ref%
{lemma_Levy_distance_difference_coefficients}}: Let $X=\sum_{k=1}^{N}\left(
d_{k}-c_{k}\right) M_{k}.$ By Lemma \ref{lemma_bound_lin_comb_Poisson_rv}, $%
\left\Vert X\right\Vert \leq b\left( 2K+3\sqrt{K}+1\right) ,$ where $%
b=\sup_{1\leq k\leq N}\left( d_{k}-c_{k}\right) $ and $K$ is the sum of the
intensities of $M_{k}.$ By Lemma \ref{lemma_Levy_distance_norm_difference},
this implies that $d_{L}\left( F_{l},F_{u}\right) \leq b\left( 2K+3\sqrt{K}%
+1\right) .$ QED.

Using Proposition \ref{lemma_Levy_distance_difference_coefficients}\textbf{,}
we can proceed to the proof of Theorem \ref%
{theorem_weak_convergence_compact_cont}. By Proposition \ref%
{theorem_weak_convergence_simple}, we know that if $N$ is fixed and $%
n\rightarrow \infty $, then 
\begin{equation*}
\sum_{i=1}^{n}l^{N}\left( X_{i,n}\right) \overset{d}{\rightarrow }%
\sum_{k=1}^{N}c_{i}^{\left( N\right) }M\left( A_{k}^{\left( N\right)
}\right) ,
\end{equation*}%
and 
\begin{equation*}
\sum_{i=1}^{n}u^{N}\left( X_{i,n}\right) \overset{d}{\rightarrow }%
\sum_{k=1}^{N}d_{i}^{\left( N\right) }M\left( A_{k}^{\left( N\right)
}\right) ,
\end{equation*}%
where $M$ is a free Poisson random measure with intensity $\lambda \left(
dx\right) .$ Let the distributions of the right-hand sides be denoted as $%
F_{l^{N}}$ and $F_{u^{N}}.$

By Corollary \ref{corollary_domination} (p. \pageref{corollary_domination}), 
$F_{l^{N}}$ is a decreasing sequence and $F_{u^{N}}$ is an increasing
sequence of distribution functions. In addition, $F_{l^{N}}\left( x\right)
\geq F_{u^{N}}\left( x\right) $ for every $N$ and $x.$ Since the sum of
intensities of variables $M\left( A_{k}^{\left( N\right) }\right) $ is less
than $\lambda \left( D\right) <\infty $ by assumption, therefore Proposition %
\ref{lemma_Levy_distance_difference_coefficients} is applicable and we can
conclude that the Levy distance between $F_{l^{N}}$ and $F_{u^{N}}$
converges to zero as $N\rightarrow \infty .$ Consequently, these two
distributions (weakly) converge to a limit distribution function as $%
N\rightarrow \infty $.

Moreover, by the definition of the integral with respect to a free Poisson
random measure, this limit equals the distribution function of $\int f\left(
x\right) M\left( dx\right) .$

In addition, by Proposition \ref{proposition_limit_point_is_between} every
limit point of the sequence of $F_{f,n}$ is between $F_{l^{N}}$ and $%
F_{u^{N}}$ for every $N,$ and therefore the sequence of $F_{f,n}$ also
converges to the distribution function of $\int f\left( x\right) M\left(
dx\right) .$ QED.

This completes the proof of Theorem \ref%
{theorem_weak_convergence_compact_cont}.

\section{Proof of Theorem 3 \label{section_proof_of_Theorem_3}}

Recall that we defined the distribution of a free order statistic in the
following way. Let $X_{1},...,X_{n}$ be freely independent self-adjoint
random variables and let $X_{i}$ have the spectral distribution functions $%
F_{i}.$ Let 
\begin{equation*}
Y_{n}\left( t\right) =\sum_{i=1}^{n}1_{\left( t,\infty \right) }\left(
X_{i}\right) .
\end{equation*}

\textbf{Definition 2}. \textit{For every real }$k\geq 0,$\textit{\ we say
that }$F_{n,k}\left( t\right) =:E\left[ 1_{\left[ 0,k\right] }\left(
Y_{n}\left( t\right) \right) \right] $\textit{\ is the distribution function
of }$k$\textit{-th\ order statistic of the sequence }$X_{1},\ldots ,X_{n},$%
\textit{\ and that it is the }$k$\textit{-th order free extremal convolution
of distribution functions }$F_{i}.$

It is straightforward to check that in the case of commutative random
variables, this definition gives the distribution function of the usual $%
\left( \left\lfloor k\right\rfloor +1\right) $-order statistic.

In the non-commutative case, we need to check that this is a consistent
definition, and that $F_{n,k}\left( t\right) $ is indeed a probability
distribution function for each $k\geq 0.$

It is easy to see that $F_{n,k}\left( t\right) $ is non-decreasing in $t.$
Indeed, let $t^{\prime }\geq t.$ Then for each $i,$ $1_{\left( t^{\prime
},\infty \right) }\left( X_{i}\right) \leq 1_{\left( t,\infty \right)
}\left( X_{i}\right) $, and therefore, $Y\left( t^{\prime }\right) \leq
Y\left( t\right) .$ It follows that $1_{\left[ 0,k\right] }\left( Y\left(
t^{\prime }\right) \right) \geq 1_{\left[ 0,k\right] }\left( Y\left(
t\right) \right) ,$ and therefore $F_{n,k}\left( t^{\prime }\right) \geq
F_{n,k}\left( t\right) .$

This function is also right-continuous in $t$. Consider a sequence $%
t_{m}\downarrow t.$ First, note that $1_{\left( t_{m},\infty \right) }\left(
X_{i}\right) \overset{d}{\rightarrow }1_{\left( t,\infty \right) }\left(
X_{i}\right) $. Second, since operators $1_{\left( t,\infty \right) }\left(
X_{i}\right) $ are freely independent for diffferent $i,$ this implies that $%
Y\left( t_{m}\right) \overset{d}{\rightarrow }Y\left( t\right) $ as $%
t_{m}\downarrow t$. Indeed, the operators $Y\left( t_{m}\right) $ and $%
Y\left( t\right) $ are uniformly bounded ($\left\Vert Y\left( t_{m}\right)
\right\Vert \leq n$ and $\left\Vert Y\left( t\right) \right\Vert \leq n$),
and the moments of the distribution of $Y\left( t_{m}\right) $ converge to
the corresponding moments of the distribution of $Y\left( t\right) .$

Third, let the spectral probability distribution functions of $Y\left(
t_{m}\right) $ and $Y\left( t\right) $ be denoted as $G_{m}\left( x\right) $
and $G\left( x\right) ,$ respectively. Then $E\left[ 1_{\left[ 0,k\right]
}\left( Y\left( t_{m}\right) \right) \right] =G_{m}\left( k\right) $ and $E%
\left[ 1_{\left[ 0,k\right] }\left( Y\left( t\right) \right) \right]
=G\left( k\right) .$ Since $G_{m}\left( k\right) \equiv F_{n,k}\left(
t_{m}\right) ,$ and $G\left( k\right) \equiv F_{n,k}\left( t\right) ,$
therefore we aim to prove that $G_{m}\left( k\right) \rightarrow G\left(
k\right) $ as $m\rightarrow \infty $ for all $k.$ The convergence $Y\left(
t_{m}\right) \overset{d}{\rightarrow }Y\left( t\right) $ means the
convergence of the moments of the spectral probability measures of operators 
$Y\left( t_{m}\right) $ and $Y\left( t\right) ,$ which implies weak
convergence of these measures because the measures have uniformly bounded
support. This implies that $G_{m}\left( k\right) \rightarrow G\left(
k\right) $ as $m\rightarrow \infty ,$ for all points $k$ at which the
probability distribution function $G\left( k\right) $ is continuous. We will
prove that, moreover, even if $G\left( x\right) $ has a jump at $x=k,$ then
the sequence $G_{m}\left( k\right) $ still converges to $G\left( k\right) .$
At this point of the argument, it is essential that $t_{m}$ converges to $t$
from above and therefore $G_{m}\left( k\right) \geq G\left( k\right) .$

Indeed, by seeking a contradiction, suppose that $G_{m}\left( k\right) $
does not converge to $G\left( k\right) .$ Then, take $\varepsilon $ such
that (i) $G_{m}\left( k\right) -G\left( k\right) >\varepsilon $ for all $m,$
and take $k^{\prime }>k$ such that (ii) $k^{\prime }$ is a point of
continuity of $G\left( x\right) ,$ and (iii) $G\left( k^{\prime }\right)
-G\left( k\right) <\varepsilon /2.$ Such $k^{\prime }$ exists because $%
G\left( x\right) $ is a spectral probability distribution function and
therefore it is right-continuous. Since $G_{m}\left( k\right) $ is
increasing, we conclude from (i), (ii), and (iii) that $G_{m}\left(
k^{\prime }\right) -G\left( k^{\prime }\right) >\varepsilon /2$ for all $m.$
But this means that $G_{m}\left( x\right) $ does not converge to $G\left(
x\right) $ at a point of continuity of $G\left( x\right) ,$ namely, at $%
k^{\prime }.$ This is a contradiction, and we conclude that $G_{m}\left(
k\right) $ converges to $G\left( k\right) $ for all $k.$ This means that $%
F_{k}\left( t\right) $ is right-continuous in $t.$

Finally, as $t\rightarrow \infty ,$ $1_{\left( t,\infty \right) }\left(
X_{i}\right) \overset{d}{\rightarrow }0.$ Therefore $Y\left( t\right) 
\overset{d}{\rightarrow }0,$ and $1_{\left[ 0,k\right] }\left( Y\left(
t\right) \right) \overset{d}{\rightarrow }I.$ Hence $F_{n,k}\left( t\right)
\rightarrow 1$ as $t\rightarrow \infty $, and we conclude that $%
F_{n,k}\left( t\right) $ is a valid distribution function.

Consider now the special case when $k=0.$ In this case $F_{n,0}\left(
t\right) $ is the dimension of the nill-space of $Y_{n}\left( t\right) ,$
which equals to the dimension of the intersection of the nill-spaces of $%
1_{\left( t,\infty \right) }\left( X_{i}\right) .$ It is easy to see that
this coincides with the definition of the free extremal convolution of the
distribution $F,$ which was introduced in \cite{benarous_voiculescu06}.

Now let us investigate the question of the limiting behavior of the
distributions $F_{n,k}\left( t\right) $ when $n\rightarrow \infty .$ The
limits are described in Theorem \ref{theorem_free_extremal_limits}.

\textbf{Proof of Theorem \ref{theorem_free_extremal_limits}:} For each $n$
we re-define:%
\begin{equation*}
Y_{n}\left( t\right) =\sum_{i=1}^{n}1_{\left( t,\infty \right) }\left( \frac{%
X_{i}-b_{n}}{a_{n}}\right) =\left\langle M_{n},1_{\left( t,\infty \right)
}\right\rangle ,
\end{equation*}%
where $M_{n}$ is the free point process associated with the triangular array 
$\left( X_{i}-b_{n}\right) /a_{n}.$

The bracket $\left\langle M_{n},1_{\left( t,\infty \right) }\right\rangle $
converges in distribution to a random variable $C_{t}$, which is a free
Poisson random variable with the intensity $\lambda \left( t\right) =-\log
G\left( t\right) .$ Then, in order to calculate the limit of $F_{n,k}\left(
t\right) $ for $n\rightarrow \infty ,$ we only need to calculate $E1_{\left[
0,k\right] }\left( C_{t}\right) ,$ that is, the distribution function of $%
C_{t}$ at $k.$ Let us denote the distribution function of $C_{t}$ as $%
H_{t}\left( x\right) ,$

For $k<0,$ we have $H_{t}\left( k\right) =0.$ For $k=0,$%
\begin{equation*}
H_{t}\left( 0\right) =\left\{ 
\begin{array}{cc}
1-\lambda \left( t\right) , & \text{ if }\lambda \left( t\right) \leq 1, \\ 
0, & \text{if }\lambda \left( t\right) >1.%
\end{array}%
\right.
\end{equation*}%
For $k>0,$ 
\begin{equation*}
H_{t}\left( k\right) =\left\{ 
\begin{array}{cc}
H_{t}\left( 0\right) , & \text{if }k<\left( 1-\sqrt{\lambda \left( t\right) }%
\right) ^{2}, \\ 
H_{t}\left( 0\right) +\int_{\left( 1-\sqrt{\lambda \left( t\right) }\right)
^{2}}^{k}p_{t}\left( \xi \right) d\xi , & \text{ if }k\in \left[ \left( 1-%
\sqrt{\lambda \left( t\right) }\right) ^{2},\left( 1+\sqrt{\lambda \left(
t\right) }\right) ^{2}\right] , \\ 
1 & \text{if }k>\left( 1+\sqrt{\lambda \left( t\right) }\right) ^{2}.%
\end{array}%
\right.
\end{equation*}%
where 
\begin{equation*}
p_{t}\left( \xi \right) =\frac{\sqrt{4\xi -\left( 1-\lambda \left( t\right)
+\xi \right) ^{2}}}{2\pi \xi }.
\end{equation*}

Then, we need to compute $F_{\left( k\right) }\left( t\right) ,$ which is $%
H_{t}\left( k\right) $ considered a function of $t$ for a fixed $k.$ Let $%
\lambda ^{-1}\left( x\right) $ denote the solution of the equation $\lambda
(t)=x.$ (That is, if $G^{\left( -1\right) }\left( x\right) $ is the
functional inversion of the limit distribution function $G\left( t\right) ,$
then $\lambda ^{-1}\left( x\right) =G^{\left( -1\right) }\left(
e^{-x}\right) .$)

Then, for $k=0$:%
\begin{equation*}
F_{\left( k\right) }\left( t\right) =\left\{ 
\begin{array}{cc}
0, & \text{ if }t\leq \lambda ^{-1}\left( 1\right) , \\ 
1-\lambda \left( t\right) , & \text{if }t>\lambda ^{-1}\left( 1\right) .%
\end{array}%
\right.
\end{equation*}%
For $k\in \left( 0,1\right) $:

\begin{equation*}
F_{\left( k\right) }\left( t\right) =\left\{ 
\begin{array}{cc}
0, & \text{if }t<\lambda ^{-1}\left( \left( 1+\sqrt{k}\right) ^{2}\right) ,
\\ 
\int_{\left( 1-\sqrt{\lambda \left( t\right) }\right) ^{2}}^{k}p_{t}\left(
\xi \right) d\xi , & \text{if }t\in \left[ \lambda ^{-1}\left( \left( 1+%
\sqrt{k}\right) ^{2}\right) ,\lambda ^{-1}\left( 1\right) \right] , \\ 
1-\lambda \left( t\right) +\int_{\left( 1-\sqrt{\lambda \left( t\right) }%
\right) ^{2}}^{k}p_{t}\left( \xi \right) d\xi , & \text{if }\left( t\in
\lambda ^{-1}\left( 1\right) ,\lambda ^{-1}\left( \left( 1-\sqrt{k}\right)
^{2}\right) \right] , \\ 
1-\lambda \left( t\right) , & \text{if }t>\lambda ^{-1}\left( \left( 1-\sqrt{%
k}\right) ^{2}\right) .%
\end{array}%
\right.
\end{equation*}

For $k\geq 1,$ we have:%
\begin{equation*}
F_{\left( k\right) }\left( t\right) =\left\{ 
\begin{array}{cc}
0, & \text{if }t<\lambda ^{-1}\left( \left( 1+\sqrt{k}\right) ^{2}\right) ,
\\ 
\int_{\left( 1-\sqrt{\lambda \left( t\right) }\right) ^{2}}^{k}p_{t}\left(
\xi \right) d\xi , & \text{if }t\in \left[ \lambda ^{-1}\left( \left( 1+%
\sqrt{k}\right) ^{2}\right) ,\lambda ^{-1}\left( 1\right) \right] , \\ 
1-\lambda \left( t\right) +\int_{\left( 1-\sqrt{\lambda \left( t\right) }%
\right) ^{2}}^{k}p_{t}\left( \xi \right) d\xi , & \text{if }\left( t\in
\lambda ^{-1}\left( 1\right) ,\lambda ^{-1}\left( \left( 1-\sqrt{k}\right)
^{2}\right) \right] , \\ 
1, & \text{if }t>\lambda ^{-1}\left( \left( 1-\sqrt{k}\right) ^{2}\right) .%
\end{array}%
\right.
\end{equation*}%
Combinining these cases, we obtain the following equation:%
\begin{equation*}
F_{\left( k\right) }\left( t\right) =\left\{ 
\begin{array}{cc}
0, & \text{if }t<\lambda ^{-1}\left( \left( 1+\sqrt{k}\right) ^{2}\right) ,
\\ 
\int_{\left( 1-\sqrt{\lambda \left( t\right) }\right) ^{2}}^{k}p_{t}\left(
\xi \right) d\xi , & \text{if }t\in \left[ \lambda ^{-1}\left( \left( 1+%
\sqrt{k}\right) ^{2}\right) ,\lambda ^{-1}\left( 1\right) \right] , \\ 
1-\lambda \left( t\right) +\int_{\left( 1-\sqrt{\lambda \left( t\right) }%
\right) ^{2}}^{k}p_{t}\left( \xi \right) d\xi , & \text{if }\left( t\in
\lambda ^{-1}\left( 1\right) ,\lambda ^{-1}\left( \left( 1-\sqrt{k}\right)
^{2}\right) \right] , \\ 
1-\lambda \left( t\right) 1_{\left[ 0,1\right) }\left( k\right) , & \text{if 
}t>\lambda ^{-1}\left( \left( 1-\sqrt{k}\right) ^{2}\right) .%
\end{array}%
\right.
\end{equation*}

QED.

\begin{example}
Distributions from the domain of attraction of Type II extremal value law
\end{example}

Consider the case of convergence to the Type II extremal value law, when the
constants $a_{n}$ and $b_{n}$ are chosen in such a way, that the limit law
is $G\left( x\right) =\exp \left( -x^{-\nu }\right) $ for $x>0.$Then we can
conclude that the limit distribution of the $k$ order statistic is given as
follows:

\begin{equation*}
F_{\left( k\right) }\left( t\right) =\left\{ 
\begin{array}{cc}
0, & \text{if }t<\left( 1+\sqrt{k}\right) ^{-2/\nu }, \\ 
\int_{\left( 1-t^{-\nu /2}\right) ^{2}}^{k}p_{t}\left( \xi \right) d\xi , & 
\text{if }t\in \left[ \left( 1+\sqrt{k}\right) ^{-2/\nu },1\right] , \\ 
1-t^{-\nu }+\int_{\left( 1-t^{-\nu /2}\right) ^{2}}^{k}p_{t}\left( \xi
\right) d\xi , & \text{if }t\in \left( 1,\left( \left( 1-\sqrt{k}\right)
^{2}\right) ^{-1/\nu }\right] , \\ 
1-t^{-\nu }1_{\left[ 0,1\right) }\left( k\right) , & \text{if }t>\left(
\left( 1-\sqrt{k}\right) ^{2}\right) ^{-1/\nu },%
\end{array}%
\right.
\end{equation*}%
where 
\begin{equation*}
p_{t}\left( \xi \right) =\frac{\sqrt{4\xi -\left( 1-t^{-\nu }+\xi \right)
^{2}}}{2\pi \xi }.
\end{equation*}

We illustrate this result for some particular values of $\nu $ and $k.$

Consider $k=0.$ Then 
\begin{equation*}
F_{\left( 0\right) }\left( t\right) =\left\{ 
\begin{array}{cc}
0, & \text{if }t<1, \\ 
1-t^{-\nu }, & \text{if }t\geq 1.%
\end{array}%
\right.
\end{equation*}%
This is the Type 2 (\textquotedblleft Pareto\textquotedblright ) limit
distribution in Definition 6.8 of \cite{benarous_voiculescu06}.

The distributions of $k$ order statistics for different values of $k$ are
illustrated in Figure 1.

\begin{center}
[Put Figure 1 here.]
\end{center}

It is interesting to note that if $k>1,$ then for all sufficiently large $t,$
$F_{\left( k\right) }\left( t\right) =1.$ This can be interpreted as saying
that the scaled $k$ order statistic is guaranteed to be less than $t_{0}$
for a suffiiciently large $t_{0}.$ In another interpretation, this result
means that for our choice of scaling parameters $a_{n}$ and $b_{n}$ and for
every $k>1,$ if $t$ is sufficiently large, then%
\begin{equation*}
\left\Vert \sum_{i=1}^{n}1_{\left( a_{n}t+b_{n},\infty \right) }\left(
X_{i}\right) \right\Vert <k
\end{equation*}%
for all large $n.$

A similar situation occurs in the classical case if the initial distribution
(i.e. the distribution of $X_{i}$) is bounded from above. In this case the
limit distribution is also bounded from above. In contrast, in the free
probability case this situation occurs even if the initial distribution is
unbounded from above. Our previous example shows that this situation occurs
even if the initial distribution has heavy tails.

\bigskip

\textbf{Acknowledgements.} \textit{The authors would like to thank two very
diligent referees for their thorough reports. G. Ben Arous woud like to
thank D. Voiculescu for useful discussions.}

\end{document}